\documentclass[12pt,a4paper]{article}
\usepackage[centertags]{amsmath}
\usepackage{amssymb}
\usepackage[utf8]{inputenc}
\usepackage{amsthm}
\usepackage{stackrel}
\usepackage{amsfonts}
\usepackage{cancel}
\usepackage{amsmath}
\usepackage[normalem]{ulem}
\usepackage{textcomp}
\usepackage[margin=3.5cm]{geometry}
\usepackage[dvips]{graphicx}
\usepackage{color}
\usepackage[dvipsnames]{xcolor}
\newdimen\AAdi%
\newbox\AAbo%
\def\AArm{\fam0 \mathrm}%
\def\AAk#1#2{\setbox\AAbo=\hbox{#2}\AAdi=\wd\AAbo\kern#1\AAdi{}}%
\def\AAr#1#2#3{\setbox\AAbo=\hbox{#2}\AAdi=\ht\AAbo\raise#1\AAdi\hbox{#3}}%
\def\BBn{{\AArm I\!N}}%
\def\BBv{{\AArm V}}%
\def\BBone{{\AArm 1\AAk{-.8}{I}I}}
\newtheorem{Theorem}{{\bf Theorem}}

\newtheorem{Proposition}{{\bf Proposition}}

\newtheorem{Lemma}{{\bf Lemma}}
\newtheorem{cor}{Corollary}

\definecolor{mypink1}{rgb}{0.858, 0.188, 0.478}

\definecolor{orangeVB}{rgb}{1.0,0.6,0.2}
\definecolor{greenVB}{rgb}{0.0, 0.5, 0.0}

\begin{document}

\hyphenation{di-stri-bu-tions} \hyphenation{e-qui-va-lent}
\hyphenation{pro-ba-bi-li-ty} \hyphenation{po-ste-rior}
\hyphenation{ac-cor-ding} \hyphenation{me-tho-do-lo-gy}
\hyphenation{mo-dels} \hyphenation{pro-per-ties}
\hyphenation{distri-buted} \hyphenation{De-ve-lo-ping}

\title{Hypotheses testing and posterior concentration rates for semi-Markov processes}
\author{V.S. Barbu\footnote{Postal address: Laboratoire de Math\'ematiques
							Rapha\"el Salem, Universit\'e de Rouen-Normandie, UMR 6085, Avenue de l'Universit\'e, BP.12,
							F76801, Saint-\'Etienne-du-Rouvray, France}, G. Gayraud\footnote{Postal address: Sorbonne University, Universit\'e de Technologie de Compi\`egne, LMAC Laboratory of Applied Mathematics of Compi\`egne - CS 60 319 - 60 203 Compi\`egne  cedex, France.}, N. Limnios$^{\dag}$, 
I. Votsi\footnote{Postal address: Laboratoire Manceau de Math\'{e}matiques, Le Mans Universit\'{e},
72000, Le Mans, France.}}

\date{} 
\maketitle

\begin{abstract}
In this paper, we adopt a nonparametric Bayesian approach and investigate the asymptotic behavior of the posterior distribution in continuous time and general state space semi-Markov processes. 
In particular, we obtain posterior concentration rates for semi-Markov kernels. For the purposes of this study, we construct robust statistical tests between Hellinger balls around semi-Markov kernels and present some specifications to particular cases, including discrete-time semi-Markov processes  and finite state space Markov processes. The objective of this paper is to provide sufficient conditions on priors and semi-Markov kernels that enable us to establish posterior concentration rates.
\end{abstract}


\section{Introduction}\label{sec:Intro}

Semi-Markov processes (SMPs) are stochastic processes that are
widely used to model real-life phenomena encountered
in seismology, biology, reliability, survival analysis,  wind energy, finance  and  other scientific fields. SMPs (\cite{LE54},\cite{SM55},\cite{T54}) generalize Markov processes in the sense that they allow the sojourn times in states to follow any distribution on $[0,+\infty)$, instead of the exponential distribution in the Markov case. Since no memoryless distributions could be considered in a semi-Markov environment, duration effects could be reproduced. The duration effect firms that
the time the semi-Markov system spends in a state influences its transition probabilities.
Particular cases of SMPs include continuous and discrete-time Markov chains and ordinary,
modified and alternating renewal processes. The foundations of the theory of SMPs were laid by Pyke (\cite{P61a}, \cite{P61b}). Since then, further significant results were obtained by
\c{C}inlar \cite{C69}, Korolyuk et al. \cite{KL05} and many others. We refer the interested reader to
Limnios and Opri\c{s}an \cite{LO01} for an approach to SMPs and their applications in reliability. For an overview in the theory on semi-Markov chains oriented toward applications in modeling and estimation see Barbu and Limnios \cite{BL08}.

Although the statistical inference of SMPs has been extensively studied from a frequentist point of view, the Bayesian literature is rather limited. Except from some specific SMP models (\cite{ECO14},\cite{EPI14}), only a few papers have considered the nonparametric Bayesian theory supporting these models (\cite{BUL07},\cite{PHE90}). Here we aim to close the aforementioned gap and follow a nonparametric Bayesian approach. The key quantity in the theory of SMPs is the semi-Markov kernel (SMK), $Q$. Our objective is to draw Bayesian inference on the Radon-Nikodym derivative of the SMK, $q$. Let us denote by $\mathcal{H}_{n}$ a trajectory of the SMP of length $n$  and by $\Pi$ the prior distribution of $q$, which in all generality, could depend on $n$, and thereafter will be denoted by $\Pi_{n}$. Given $\mathcal{H}_{n}$ and $\Pi_{n},$ the knowledge on $q$ is updated by the posterior distribution, that is denoted by $\Pi_{n}^{\mathcal{H}_{n}}(\cdot)=\Pi_{n}(\cdot|{\mathcal{H}_{n}}).$  We shall stick to the last notation throughout the paper and further denote by
$q_{0}$ the derivative of the ``true'' SMK, $Q_{0}$, which is the SMK that generated $\mathcal{H}_{n}$. The main topic of the article is the study of the asymptotic behaviour of $\Pi_{n}^{\mathcal{H}_{n}}$ in a neighbourhood of $Q_0.$

Most of the known results in the asymptotic behaviour of posterior distributions in infinite-dimensional models address issues of the posterior consistency and posterior concentration around the true distribution. In a nonparametric context, when the observations are i.i.d., such results were first derived in \cite{GGD00} and \cite{SW01} with a variety of examples. Beyond the i.i.d. setup, the asymptotic  behaviour of the posterior has been studied in the context of independent nonidentically distributed observations (\cite{AA03}, \cite{AGR13}, \cite{CGR04}, \cite{GDV07}, \cite{GR06}, \cite{GGR99}). 




One of the most natural extensions of the i.i.d. structure is a Markov process, where only the immediate past matters. Although, given the present, the future will not further depend on the past, the dependence propagates and may reasonably capture the dependence structure of the observations. Ghosal and van Der Vaart \cite{GDV07} studied the asymptotic behaviour of posterior distributions to several classes of non-i.i.d. models including Markov chains. For their purposes the authors used previous results on the existence of statistical tests (\cite{B83a}, \cite{LC86}, \cite{LC75}, \cite{LC73})
between two Hellinger balls for a given class of models. We refer the interested reader to \cite{B2013} for improved results about the existence of such tests for the relevant estimation
problems. Tang and Ghosal \cite{TG07} extended Schwartz's theory of posterior consistency to ergodic Markov processes and applied it in the context of a Dirichlet mixture model for transition densities.  More recently, Gassiat and Rousseau \cite{GR13} studied the posterior distribution in hidden Markov chains where both the observational and the state spaces are general. For nonparametric Bayesian estimation of conditional distributions, Pati et al. \cite{PAT13} provided sufficient conditions on the prior under which the weak and various types of strong posterior consistency could be obtained.

For reviews on posterior consistency as well as posterior concentration in infinite dimensions, the interested reader can refer to Wasserman \cite{W98}, Ghosh and Ramamoorthi \cite{GR97} and Ghosal et al. \cite{GDV17}. 

This paper aims to extend previous results by studying the convergence of the posterior distribution
of $q$ for SMPs. Specifically, we generalize and extend previous results on discrete-time Markov processes in finite state space \cite{GDV07} to continuous-time SMPs in general state space.

In order to apply the general theory to the semi-Markov framework,  we demonstrate the existence of the relevant statistical tests. To this purpose, we extend the hypotheses testing results  for Marov chains developed by  Birg\'{e} \cite{B83a} to continuous-time general state space SMPs. Such tests can also be used to distinguish Markov from semi-Markov models and decide which model could better describe the data, which is  a crucial subject in real-world applications. 

Very few researchers considered hypotheses testing problem in a semi-Markov context. Bath and Deshpande \cite{bath} developed a nonparametric test for testing Markov against semi-Markov processes. Banerjee and Bhattacharyya \cite{baner} considered a two-state SMP and proposed parametric tests for the equality of the sojourn time distributions, under the assumption that these distributions are absolutely continuous and belong to the Exponential family. Also in a parametric context, Malinovskii \cite{MAL92} considered that the probability distribution of an SMP depends on a real-valued parameter $\vartheta>0$ and studied the simple hypothesis $H_0: \vartheta=0$ against $H_{1}: \vartheta=hT^{-1/2}$, $0<h\leq{c}$ (the SMP is observed up to time $T$). Chang et al. (\cite{chang1999}, \cite{chang2001}) considered hypotheses testing problems for semi-Markov counting processes, in a survival analysis context. Tsai  \cite{tsai} proposed a rank test based on semi-Markov processes in order to test whether a pair of observation $(X, Y)$ has the same distribution as $(Y, X),$ i.e., $X,$ $Y$ exchangeable.  To the best of our knowledge, the present research is the first one that considers general robust hypotheses testing problems for SMPs in a nonparametric context. 

We focus on SMPs since they are much more general and better adapted to applications than the Markov
processes.  In real-world systems, the state space of the under study processes could be $\{0,1\}^{\mathbb{N}},$ (e.g., communication systems), where $\mathbb{N}$ is the set of nonnegative integers, 
or $[0, \infty)$ (e.g., fatigue crack growth modelling). This is the reason why we concentrate on  general SMPs. On the other side, since in physical and biological applications time is usually considered to be continuous,  discrete-time processes are not always appropriate for describing such phenomena. In such situations continuous-time processes are often more suitable than the discrete-time ones. 
Therefore we focus our discussion on the continuous-time case rather than the discrete-time case. Nonetheless, note  that our results on the robust tests are very general
and could also be applied to the discrete-time case, with the corresponding modifications. 

The organization of the paper is as follows. In Section $2$ the notation and preliminaries
of semi-Markov processes are presented; the objectives of our paper are also presented. Section $3$ describes the hypotheses testing for the processes under study 
and some particular cases. Section $4$ discusses the derivation of the posterior concentration rate and the relative hypotheses. Finally, in Section $5$, we give a detailed description of the proofs and some technical lemmas.


\section{The semi-Markov framework and objectives}\label{sec:semi-M}



\subsection{Semi-Markov processes}

We consider $(E,\mathcal{E})$ a measurable space and an $(E,\mathcal{E})-$valued semi-Markov process $\mathbf{Z}:=(Z_t)_{t \in \mathbb{R}^{+}}$ defined on a complete probability space $(\Omega,\mathcal{F},\mathbb{P})$. The semi-Markov process $\mathbf{Z}$ corresponding to the 
Markov renewal process (MRP)
$\mathbf{(J,S)}:=(J_{n},S_{n})_{n\in{\mathbb{N}}}$, is defined by
\begin{equation*}
Z_{t}:=J_{N(t)}, \text{ } \text{ } \text{ } t\in{\mathbb{R}^{+}},
\end{equation*}
where
 $0\leq{S_{0}}\leq\ldots\leq{S_{n}}\leq\ldots$ are the successive 
$\mathbb{R}^{+}$-valued jump times of $\mathbf{Z}$, $(J_n)_{n\geq  0}$ denotes the successive visited states at these jump times (henceforth called \textit{the embedded Markov chain (EMC)}) and
\[N(t)=\left\{\begin{array} {ll}
              0, &
              \mbox{if} \quad S_{1}-S_{0}>t,\\
              \sup\{n\in\mathbb{N^{*}}: S_{n}\leq{t}\},  &
              \mbox{if} \quad S_{1}-S_{0}\leq{t}.
              \end{array}
         \right.
         \]
$S_0$ may be viewed as the first non-negative time at which  a jump is observed.
In what follows,  the EMC and MRP are considered to be homogeneous with respect to $n\in{\mathbb{N}}$. It is worth noticing that the MRP $\mathbf{(J,S)}$
satisfies the following Markov property, i.e., for any $n\in{\mathbb{N}}$,
any $t\in{\mathbb{R}^{+}}$ and any $B\in{\mathcal{E}}$:
\begin{equation*}
\mathbb{P}(J_{n+1}\in{B},S_{n+1}-S_{n}\leq{t}|J_{0},\ldots,J_{n},S_{0},\ldots,S_{n})\\
\stackrel{\mbox{{\rm a.s.}}}{=}\mathbb{P}(J_{n+1}\in{B},S_{n+1}-S_{n}\leq{t}|J_{n}).
\end{equation*}


In the semi-Markov framework, of central importance is the semi-Markov kernel (SMK) defined as follows:
\begin{equation*}
Q_x(B,t):= \mathbb{P}(J_{n+1}\in{B},S_{n+1}-S_{n}\leq{t}|J_{n}=x),\; x\in{E}, \; t\in{\mathbb{R}^{+}}, \; B\in{\mathcal{E}}\\
\end{equation*}
Since we suppose that the distribution of $\mathbf{Z}$ is unknown, we
focus our interest on the semi-Markov kernel. In particular the stochastic behavior of the SMP $\textbf{Z}$ is determined 
completely by its SMK and its initial distribution.

Let us denote the $n-$step transition kernel  of the  EMC $(J_{n})_{n\in{\mathbb{N}}}$ by 

\begin{equation}\label{trans-eq-J}
P^{(n)}(x,B):= \mathbb{P}(J_{n}\in{B}|J_{0}=x),\; x\in{E}, \; B \in {\cal E},
\end{equation}
and the (one-step) transition kernel by $P(x,B)=Q_x(B,\infty).$  


It is worth mentioning that
$$Q_x(B,t)= \int_B P(x,dy) \mathbb{P}(S_{n+1}- S_n \leq t |J_n=x, J_{n+1}=y), \forall t \in \mathbb{R}^+, \; \forall B \in {\cal E}.$$

The following assumptions have to be considered in the sequel.

\begin{enumerate}
  \item [\textbf{A1}] The embedded Markov chain $(J_{n})_{n\in\mathbb{N}}$ is ergodic with stationary probability measure $\boldsymbol{\rho}$ (that is
$\boldsymbol{\rho}{P}=\boldsymbol{\rho}$, with $P$ the transition kernel of $\mathbf{J}$ and
$\boldsymbol{\rho}(E)=1$).
  \item [\textbf{A2}] The mean sojourn times $m(x)=\int_{0}^{\infty}\mathbb{P}(S_{1}-S_{0}>t \mid J_0 = x)dt$
satisfies
\begin{equation*}
    \int_{E}^{}\boldsymbol{\rho}(dx)m(x)<{\infty}.
\end{equation*}

\item [\textbf{A3}]  
\begin{equation*}
\mathbb{P}(S_{n+1}- S_n \leq t |J_n=x, J_{n+1}=y)\neq
\textbf{1}_{\mathbb{R}^{+}}(t),
\forall n \in \mathbb{N}, \; \forall t \in \mathbb{R}^+, \; \forall x, y \in E.
\end{equation*}

\end{enumerate}
Note that  \textbf{A2}  and \textbf{A3} ensure that for all non negative $t$ and $B \in {\cal E}$, $\mathbb{P} (Z_t \in B)$ is always well-defined and non-zero. However the conditional probability in  Assumption \textbf{A3} may be defined as any Dirac measure on positive real numbers.

Denote also by $\mathbb{B}^{+}$ the Borelian $\sigma-$algebra on $\mathbb{R}^{+}$. We suppose that for any $x \in E$, the SMK starting from $x$ is absolutely continuous with
respect to (w.r.t.) $\nu$, a $\sigma-$finite measure 
$(E\times\mathbb{R}^{+},\mathcal{E}\otimes{\mathbb{B}^{+}})$ and denote by $q_x(\cdot,\cdot)$ its Radon-Nikodym (RN) derivative, i.e.,
$Q_x(dy,dt)=q_x(y,t)d{\nu}(y,t)$. For $n\geq{1}$, let $X_{n}:=S_{n}-S_{n-1}$ be the successive sojourn times of
$\mathbf{Z}$ and ${0}\leq{X_{0}}=S_{0}$. On $\mathcal{E}\otimes{\mathbb{B}^{+}}$, we further define  
the measure $\widetilde{\boldsymbol{\rho}}$  as the distribution of $\mathbf{(J,X)}:=(J_{n},X_{n})_{n\in{\mathbb{N}}}$, where 
\begin{equation} \label{eq:carac}
\widetilde{\boldsymbol{\rho}} (A,{\Gamma})=\int_{E}{\boldsymbol{\rho}}(dx)Q_x(A,{\Gamma}), \; \forall A \in \mathcal{E}, \; \forall \Gamma \in \mathbb{B}^+. 
\end{equation}

\begin{Proposition}\label{prop:caracterisation}
The measure  $\widetilde{\boldsymbol{\rho}}$ defined in  \eqref{eq:carac} is the stationary distribution 
of $(J_{n},X_{n})_{n\in{\mathbb{N}}}$. 
\end{Proposition}

Since we are interested in obtaining asymptotic results, without loss of generality  we consider as initial distribution of the process $\mathbf{(J,X)}$ its stationary distribution, $\widetilde{\boldsymbol{\rho}}$. To avoid complicated notation, we will also use $\widetilde{\boldsymbol{\rho}}$ to denote the density w.r.t. $\nu.$ 

In the sequel, the hypotheses \textbf{A1},  \textbf{A2}  and \textbf{A3} are considered to hold true.

\subsection{Objectives}


Recall that we have denoted  by $Q_{0}$ the true semi-Markov kernel and by 
${q}_{0}$ its RN derivative w.r.t. $\nu,$ cf. Section $2.$ We suppose that $q_0$ belongs to a certain set of semi-Markov kernel densities $\mathcal{Q}$ defined by
$$\mathcal{Q}=\{q=q_x(y,t): x,y\in{E}, t\in{\mathbb{R}^{+}}\},$$
which is equipped with a metric $d$ that will be defined 
in the sequel. Next consider $\epsilon$-neighborhoods around $q_0$ in $\mathcal{Q}$ w.r.t. $d$, that is
 \begin{equation*}
B_{d}(q_{0},\epsilon)=\Big{\{}q\in{\mathcal{Q}}:d({q}_{0},{q})\leq{\epsilon}\Big{\}}.
\end{equation*}
To allow some flexibility, it is quite common to deal with  $\mathcal{Q}_n$, a subset of $\mathcal{Q}$, that may depend on $n,$  such that the prior distribution $\Pi_n$ on ${\cal Q}$ assigns most of its mass on $\mathcal{Q}_n$ (see Assumption \textbf{H4} below). An $\epsilon$-neighborhood around $q_0$ in $\mathcal{Q}_n$ w.r.t. $d$ will be denoted by $B_{d,n}(q_{0},\epsilon)$.

\bigskip
As noted by Birg\'e \cite{B83a} in the setting of Markov chains,
there exists a priori no ``natural" distance $d$ between two semi-Markov kernel densities. Nevertheless, a natural distance could be defined between two probability distributions
$Q_{x;1}$ and $Q_{x;2}$ dominated by $\nu$ and corresponding to the same initial state $J_0=x \in E$. Indeed, if we further denote by ${q}_{x;1}$ and ${q}_{x;2}$ their respective RN derivatives, and following the lines of Birg\'e \cite{B83a}, $d$ could be defined in two steps. First by
considering the squared Hellinger distance between ${Q}_{x;1}$ and ${Q}_{x;2}$, i.e.,
\begin{eqnarray}\label{Hellinger_nu}
h_{\nu}^{2}(Q_{x;1},Q_{x;2})&=&\frac{1}{2}\int_{E \times{\mathbb{R}^{+}}}
{\Bigl{(}\sqrt{q_{x;1}(y,t)}-\sqrt{q_{x;2}(y,t)}\Bigr{)}^{2}d{\nu}(y,t)},
 \label{Hell-dist}
\end{eqnarray}
and second, given a measure on $\mathcal{E}$, say  $\mu$, by  setting a semi-distance $d_\mu$ between $q_{1}$ and $q_{2}$,
\begin{eqnarray}
d_{\mu}^{2}(q_{1},q_{2})& =& \int_{E}^{}h_{\nu}^{2}(Q_{x;1},Q_{x;2})d{\mu}(x).  \label{def:distance-kernel}
\end{eqnarray}

Given a sample path of the SMP for a given number of jumps $n\in{\mathbb{N}^{*}}$,
\begin{equation*}
\mathcal{H}_{n}=\{J_{0},J_{1},\ldots,J_{n},S_{0},S_{1},\ldots,S_{n}\},
\end{equation*}
we adopt a Bayesian point of view by considering a prior distribution $\Pi_n$ on  $\mathcal{Q}.$ We aim to establish how fast the posterior distribution shrinks, in terms of $d$, the ``true'' semi-Markov kernel density, $q_{0}$. The precise definition of $d$ will be given after the statement of Assumption \textbf{H1}, where the measure $\mu$ is fixed. More precisely,
our objective is to find the minimal positive sequence $\epsilon_n$ tending to zero as $n$ goes to infinity, such that
under some assumptions on both $\mathcal{Q}$ and $\Pi_n$
\begin{eqnarray*}
\Pi_n^{\mathcal{H}_{n}} \left( B^{\complement}_d({q}_{0},\epsilon_{n})    \right)
\stackrel{L_1(\mathbb{P}^{(n)}_0)}{\longrightarrow} 0 \text{ } \text{as} \text{ } {n\rightarrow 0},
\end{eqnarray*}
where $B^{\complement}_d$ denotes the complementary of $B_d$ in $\mathcal{Q}$ and $\mathbb{P}^{(n)}_0$ refers to the ``true" distribution of ${\cal H}_n.$

Let us denote by  $\mathbb{P}_q^{(n)}$ the distribution of $\mathcal{H}_{n}$, when the density of the SMK is $q$. We further denote by $\mathbb{E}_q^{(n)}$ the expectation and by $\BBv^{(n)}_{q}$ the variance w.r.t. $\mathbb{P}_q^{(n)}$, respectively.  Every  quantity (distribution, SMK, expectation, variance,$\ldots$) with an index $0$ refers to the corresponding ``true'' quantity.


\section{Hypotheses testing for semi-Markov processes}\label{sec:Bayes}


\subsection{Robust tests} \label{subsec:robust}

One of the key ingredients needed to obtain posterior concentration rates is the construction of corresponding robust hypotheses tests.  For a variety of models, depending on the semi-metric $d$, some tests with exponential power do exist. For instance, in the case of density or conditional density estimation, Hellinger or $L_{1}$ tests have been introduced in \cite{B83b}. Other examples of tests could be found in \cite{GDV07} and in \cite{RCL12}. However, to the best of our knowledge, no such tests exist for semi-Markov processes. Therefore it is of paramount importance to build test procedures with exponentially small errors in the semi-Markov context. Thus in the sequel we will be interested in the following testing procedure
\begin{eqnarray}\label{eq_HypoTest}
H_{0}: q_{0} \;\text{against} \; H_{1}: q \in B_{d_{\eta*},n}(q_{1},\xi\epsilon), \; \text{{\rm with}} \;d_{\nu*}(q_{0},q_{1})\geq{\epsilon},
\end{eqnarray}
for some  $\xi \in (0,1).$



In order to derive posterior concentration rates for SMK densities, one more assumption is required.
 
\begin{itemize}
\item \textbf{H1}: 
There exist two measures $\nu^{*}$ and $\eta^{*}$ on
$\mathcal{E}$ and two positive integers $k,l$ such that for any ${x}\in{E}$,
\begin{equation*}
{\frac{1}{k}\sum_{u=1}^{k}{P^{(u)}{(x,\cdot)}}}  \geq \nu^{*}(\cdot) \quad \quad \text{and} \quad \quad
{P^{(l)}(x,\cdot)\leq \eta^{*}(\cdot)},
\end{equation*}
\end{itemize}
where $P^{(\cdot)}$ is defined in \eqref{trans-eq-J}. Note that  \textbf{H1} implies the following inequalities which 
serve to prove Proposition \ref{Prop:test}: 
\begin{equation*}
\forall m\in{\mathbb{N}}, \quad \quad \frac{1}{k}\sum_{u=1}^{k}{P^{(u+m)}{(x,\cdot)}}  \geq \nu^{*}(\cdot) \quad \quad \text{and} \quad \quad
P^{(l+m)}(x,\cdot)\leq \eta^{*}(\cdot) .
\end{equation*}

\begin{Proposition}\label{Prop:test}
Under Hypothesis \textbf{H1},  for any $n\in\mathbb{N}^{*}$,
there exist
universal positive constants $\xi \in (0,1),$ $K$  and ${\tilde K}$ such that for any
$\epsilon>0$ and any $q_{1}\in\mathcal{Q}_n$ such that 
$d_{\nu*}(q_{1},q_{0})>\epsilon$,
there exists a test $\psi_{1} (\mathcal{H}_n)$ satisfying
\begin{eqnarray}\label{test:1-2-error}
\mathbb{E}_{0}^{(n)}[\psi_{1} (\mathcal{H}_n)]\leq{e^{-Kn\epsilon^{2}
}} \text{ } \text{ } \text{ } \text{{\rm and}} \text{ } \text{ }\text{ }
\sup_{q\in\mathcal{Q}_n : d_{\eta*}(q_{1},q)<\epsilon\xi}
\mathbb{E}_q^{(n)}[1-\psi_{1} (\mathcal{H}_n)]\leq{e^{-{\tilde K}n\epsilon^{2}}}.
\end{eqnarray}
\end{Proposition}
\noindent

%

The next corollary generalizes Proposition \ref{Prop:test} to any $q_1 \in {\cal Q}_n$ which is $\epsilon-$distant from $q_0$ w.r.t. $d_{\nu*}$. It requires an additional assumption (see hereafter \textbf{H2}) to control the complexity of ${\tilde {\cal Q}_n}\subseteq{{\cal Q}_n}$. This assumption is 
based on the minimum number of $d_{\nu*}$-balls of radius ${\tilde {\epsilon}}$ needed to cover ${\tilde {\cal Q}_n}$, which is denoted by ${N}({\tilde {\epsilon}},{\tilde {\cal Q}_n},d_{\nu*})$.

Note that the case where the null hypothesis is composite could also be considered; the first type error in \eqref{test:1-2-error} would be written similarly to the second type error, with straightforward modifications.  

\begin{cor}\label{cor:test-global}
Under Hypothesis \textbf{H1}, assume that for a sequence $\epsilon_n$ of positive numbers such that $\displaystyle{\lim_{n \rightarrow +\infty}} \epsilon_n = 0$ and $\displaystyle{\lim_{n \rightarrow +\infty}} n \epsilon_n^2=0$, the following assumption holds true.

\begin{itemize}
\item \textbf{H2} 
%
For  ${\xi}$ in $(0,1),$ 
$$\sup_{\epsilon{>}\epsilon_{n}}\log{N}\bigl{(}{\epsilon\xi},B_{d_{\nu*},n}(q_0,\epsilon),d_{\eta*}\bigr{)}\leq{n\epsilon_n^{2}}.$$ 
\end{itemize} 
Then, 
there exists a test $\psi(\mathcal{H}_n)$ satisfying 
\begin{align*}
& \mathbb{E}_{0}^{(n)}[\psi(\mathcal{H}_n)]\leq{e^{-Kn\epsilon^{2}_{n}M^{2}}}\\
\text{ } \text{ } \text{{\rm and} } \text{ } \text{ } & \nonumber \\ 
& \sup_{q \in \mathcal{Q}_n:  d_{\nu*}(q_{0},q)>\epsilon_{n} M}
\mathbb{E}_q^{(n)}[1-\psi(\mathcal{H}_n)]\leq{e^{-{\tilde K}n\epsilon^{2}_{n}M^{2}}}.
\end{align*}
\end{cor}

\subsection{Particular cases}

In this paper the results are rather generic in the sense that 
they refer to continuous-time and general state space SMPs. In the sequel, we focus on some particular cases
that could be of special interest either from an applicative point of view, or as a starting point for further research.

First, note that the state space is considered to be finite in most of the applicative articles.
%
Second, we would like to stress out that in some applications the state space is intrinsically continuous, due to the fact that the scale of the measures is continuous.

\subsubsection{Discrete-time SMPs}

\begin{itemize}

\item \textsc{General state space}

Let us first denote by
\begin{equation*}
q_{x}(y,k)=\mathbb{P}(J_{n+1}=y,X_{n+1}=k|J_{n}=x),
\end{equation*}
the RN derivative of the SMK. 
Then for any $k\in{\mathbb{N}}$ and any $B\in{\mathcal{E}}$, 
the respective cumulative semi-Markov kernel is given by
\begin{equation*}
Q_{x}(B,k)=\mathbb{P}(J_{n+1}\in{B},X_{n+1}\leq{k}|J_{n}=x).
\end{equation*}

It should be noted that in this case $\nu$ in \eqref{Hellinger_nu} is the product measure between  a finite-measure $\mu$ on $(E,\mathcal{E})$ used in (\ref{def:distance-kernel}) and the counting measure on $\mathbb N.$ Thus in this framework, the squared Hellinger distance becomes
\begin{eqnarray*}
h_{\mu}^{2}(Q_{x;1},Q_{x;2})&=&\frac{1}{2}\sum_{k\in{\mathbb{N}}}^{}\int_{E}
{\Bigl{(}\sqrt{q_{x;1}(y,k)}-\sqrt{q_{x;2}(y,k)}\Bigr{)}^{2}d{\mu}(y)},
\end{eqnarray*}
while the semi-distance $d_\mu$ between $q_{1}$ and $q_{2}$ is 
given in Equation (\ref{def:distance-kernel}). 


\item \textsc{Finite state space}

For any $k\in{\mathbb{N}}$ and any $y\in{E}$, we define by
\begin{equation}\label{eq_SMK_FiStDT}
q_{x}(y,k)=\mathbb{P}(J_{n+1}=y,X_{n+1}=k|J_{n}=x),
\end{equation}
the semi-Markov kernel and by
\begin{equation*}
Q_{x}(y,k)=\mathbb{P}(J_{n+1}=y,X_{n+1}\leq{k}|J_{n}=x)
\end{equation*}
the cumulative semi-Markov kernel, respectively. 

Since in this framework $\mu$ is the counting measure on  $(E,\mathcal{E}),$ the squared Hellinger distance becomes
\begin{eqnarray}\label{eq_Hellinger_FiStDT}
h^{2}(Q_{x;1},Q_{x;2})&=&\frac{1}{2}\sum_{k\in{\mathbb{N}}}^{}\sum_{y\in{E}}^{}
\Bigl{(}\sqrt{q_{x;1}(y,k)}-\sqrt{q_{x;2}(y,k)}\Bigr{)}^{2},
\end{eqnarray}
and the semi-distance $d$ between $q_{1}$ and $q_{2}$ is given by
\begin{eqnarray}
d^{2}(q_{1},q_{2})& =& \sum_{x\in{E}}^{}h^{2}(Q_{x;1},Q_{x;2}).  
\label{semi_dist}
\end{eqnarray}

\end{itemize}

\subsubsection{Continuous-time SMPs}

\begin{itemize}

\item \textsc{Finite state space}

Let us first denote by
\begin{equation}\label{eq_SMK_FiStCT}
Q_{x}(y,t)=\mathbb{P}(J_{n+1}=y,X_{n+1}\leq{t}|J_{n}=x)
\end{equation}
the semi-Markov kernel, for any $y\in{E}$ and any $t \in \mathbb{R}^{+}$.

In this context, the squared Hellinger distance becomes
\begin{eqnarray*}\label{eq_Hellinger_FiStCT}
h_{\nu_1}^{2}(Q_{x;1},Q_{x;2})&=&\frac{1}{2}\sum_{y\in{E}
}^{}\int_{\mathbb{R}^{+}}{\Bigl{(}\sqrt{q_{x;1}(y,t)}-\sqrt{q_{x;2}(y,t)}\Bigr{)}^{2}d{\nu_1}(t)},
\end{eqnarray*}
where $\nu_1$ is the marginal on $(\mathbb{R}^{+},\mathbb{B}^{+})$ of the measure $\nu$ defined on $E \times \mathbb{R}^{+},$
and the semi-distance $d$ between $q_{1}$ and $q_{2}$ is defined as in Eq. (\ref{semi_dist}).


%

\subsection{Specification to the Markov case}

Note that the previously obtained results on robust tests for SMPs could be adapted to the particular case of Markov processes. These tests are of great interest and
could be used for real-life applications. In particular, they  enable us to decide if an observed dataset would be better described by a Markov (null hypothesis) or a semi-Markov process (alternative hypothesis). More precisely suppose we are interested in the following testing problem  
\begin{eqnarray*}
&&\tilde{H}_0 : Q_0 \text{ Markov kernel }  \text{ vs } \\ 
&&\tilde{H}_1 : Q_1 \text{ semi-Markov kernel } \epsilon \text{ distant from } Q_0 \text{ w.r.t.
some pseudo-metric.} 
\end{eqnarray*}
Note that $\tilde{H}_1$ could be extended to any $\xi \epsilon-$ball around $Q_1$ with $\xi \in ]0,1[$. 

%

In this section, we are going to explain how the hypothesis testing problem $\tilde{H}_0$ versus $\tilde{H}_1$ can directly be
handled from solving the hypothesis problem ${H}_0$ versus ${H}_1$ stated in \eqref{eq_HypoTest}.

First, for the discrete-time and finite state space case, assume that we have a Markov process with Markov transition matrix $\widetilde{p}=(\widetilde{p}_{xy})_{x, y \in E},$  $\widetilde{p}_{xx}\neq 1$ for all states $x \in E.$ 


Note that a Markov process could be represented as a semi-Markov process with semi-Markov kernel given in \eqref{eq_SMK_FiStDT} and expressed as 
\begin{eqnarray*}
q_{x; 0}(y,k)
&=&\label{qMasSM}
\left\{ \begin{array}{ll}
\widetilde{p}_{xy}\,(\widetilde{p}_{xx})^{k-1}, & \textrm{if } x\neq y \textrm{ and } k \in \mathbb N^*,\\
0, & \textrm{otherwise.}
\end{array} \right.\\
\end{eqnarray*}
Consequently,  we can define the corresponding  squared Hellinger distance as in \eqref{eq_Hellinger_FiStDT} and construct the corresponding testing procedure. \\

Second, for the continuous-time and finite state space case, consider a regular jump Markov process with continuous transition semigroup $\widetilde{P} = \left(\widetilde{P}(t)\right)_{t \in \mathbb{R}^+}$ and infinitesimal generator matrix $A=(a_{xy})_{x, y \in E}.$ 

%
In this context, we can represent the Markov process as a semi-Markov process with semi-Markov kernel given in \eqref{eq_SMK_FiStCT} and expressed as 
\begin{eqnarray*}
Q_{x;0}(y, t)&=&\label{qMasSM_cont}
\left\{ \begin{array}{ll}
\frac{a_{xy}}{a_x} (1 - \exp(-a_x t)), & \textrm{if } x\neq y \textrm{ and } t \in \mathbb R^+,\\
0, & \textrm{otherwise,}
\end{array} \right.
\end{eqnarray*}
where $a_x := -a_{xx} < \infty, x \in E.$
%

Note that one can also consider the case where the null hypothesis is composite or the case where the alternative hypothesis is simple, with straightforward modifications.

%

%


 
\section{Posterior concentration rates for semi-Markov kernels}\label{sec:resu-st}


In this part, we present the key assumptions and state our main result.

First note that the likelihood function of the sample path $\mathcal{H}_{n}$
evaluated at $q \in {\cal Q}$  is given by
\begin{eqnarray*}
\mathcal{L}_{n}(q) &=&\widetilde{\boldsymbol{\rho}}(J_{0},S_{0}){\prod_{\ell=1}^{n}q_{J_{\ell-1}}(J_{\ell},X_{\ell})}.
\end{eqnarray*}

Let us introduce the tools that play a central role in asymptotic Bayesian nonparametrics: the Kullback-Liebler (KL) divergence between any two distributions $\mathbb{P}^{(n)}_{q_{1}}$ and $\mathbb{P}^{(n)}_{q_{2}}$ and the centered second moment of the integrand of the corresponding KL divergence, which are defined by 
\begin{eqnarray*}
K(\mathbb{P}^{(n)}_{q_{1}},\mathbb{P}^{(n)}_{q_{2}})& :=& \mathbb{E}_{0}^{(n)}\Big{[}\log\frac{\widetilde{\boldsymbol{\rho}}_{1}
(J_{0},S_{0})}{\widetilde{\boldsymbol{\rho}}_{2}(J_{0},S_{0})}\prod_{l=1}^{n}\frac{q_{J_{l-1};1}
(J_{l},X_{l})}{q_{J_{l-1};2}(J_{l},X_{l})}\Big{]},\\
V_0(\mathbb{P}^{(n)}_{q_{1}},\mathbb{P}^{(n)}_{q_{2}})& := &
\mathbb{V}_{0}^{(n)}\Big{[}\log\frac{\widetilde{\boldsymbol{\rho}}_{1}
(J_{0},S_{0})}{\widetilde{\boldsymbol{\rho}}_{2}(J_{0},S_{0})}\prod_{l=1}^{n}\frac{q_{J_{l-1};1}
(J_{l},X_{l})}{q_{J_{l-1};2}(J_{l},X_{l})}\Big{]},
\end{eqnarray*}
where $\mathbb{E}_{0}^{(n)}$ and  $\mathbb{V}_{0}^{(n)}$ denote respectively the expectation and the variance w.r.t. $\mathbb{P}^{(n)}_{0}$.

Then, consider the subspace of $\mathcal{Q}$, ${\cal U}(q_{0},\epsilon)$, which   
represents the following Kullback-Liebler $\epsilon$-neighborhood of $\mathbb{P}^{(n)}_{0}$, that is, for positive $\epsilon$, 
\begin{equation*}
{\cal U}(q_{0},\epsilon)=\Big{\{}q \in\mathcal{Q} :
K(\mathbb{P}^{(n)}_{0},\mathbb{P}_q^{(n)})\leq{n\epsilon^{2}},V_0(\mathbb{P}^{(n)}_{0},
\mathbb{P}^{(n)}_q)\leq{n\epsilon^{2}}\Big{\}}.
\end{equation*}
\medskip

It is worth mentioning that although $\widetilde{\boldsymbol{\rho}}$ is not of primary interest, since it is unknown it should require a prior. But since any prior on $\widetilde{\boldsymbol{\rho}}$  that is independent of the prior on $q$ would disappear upon marginalization of the posterior of $(\widetilde{\boldsymbol{\rho}}, q)$  relatively to $\widetilde{\boldsymbol{\rho}},$ in the sequel it will be dropped. Thus, it suffices to consider only a prior distribution on $q.$

Let us now state the main result. We recall that $\Pi_{n}$ denotes a prior distribution on $\mathcal{Q}.$

\begin{Theorem}\label{THM}
Assume that \rm \textbf{H1} holds true and suppose that for a sequence of positive numbers $\epsilon_n$  such that $\displaystyle{\lim_{n \rightarrow +\infty}} \epsilon_n = 0$, $\displaystyle{\lim_{n \rightarrow +\infty}} n \epsilon_n^2=0$, \textbf{H2} and \textbf{H3-H4} defined hereafter, hold true.


\begin{itemize}

\bigskip

\item \textbf{H3}   $\exists \; c>0$, 
$\Pi_n\bigl{(}{\cal U}(q_{0},\epsilon_{n})\bigr{)}>{e^{-cn\epsilon_{n}^{2}}}$,
\bigskip
\item \textbf{H4}  $\mathcal{Q}_n \subset \mathcal{Q}$ is such that
$\Pi_{n}\bigl{(}\mathcal{Q}_n^{\complement}\bigr{)}\leq{e^{-2n(c+1)\epsilon^{2}_{n}}}.$
%
\end{itemize}

\bigskip


Then for $M$ large enough,

\begin{eqnarray}\label{eq:resu}
\Pi_{n}^{\mathcal{H}_{n}}\bigl{(}B^{\complement}_{d_{\nu*}}(q_{0},\epsilon_{n}M)\bigr{)} \stackrel{L_1 (\mathbb{P}_0^{(n)})}{\longrightarrow} 0,\quad \mbox{{\rm as}} \quad n \rightarrow \infty.
\end{eqnarray}

\end{Theorem}

Some comments on the result of Theorem \ref{THM} as well as the hypotheses we deal with:

\begin{itemize}

\item Under   {\rm \textbf{H1}}, Theorem \ref{THM} guarantees that, for both a particular set of semi-Markov kernels ${\mathcal Q}$  containing some subset 
${\mathcal Q}_n$ such that {\rm \textbf{H2}} holds true for a sequence of positive numbers  $\epsilon_n$  and a prior distribution  $\Pi_n$ on 
${\mathcal Q}$ satisfying  assumptions {\rm \textbf{H3}-\textbf{H4}} with $\epsilon_n$,  the posterior distribution shrinks towards $q_0 \in {\cal Q}$ at a rate  proportional to $\epsilon_n$. 

\item Assumption  {\rm \textbf{H3}} is classical in Bayesian Nonparametrics; it states that the prior distribution puts enough mass around 
KL  neighborhoods of $q_0$.

\item  As mentioned  in Section \ref{subsec:robust}, ${\cal Q}_n$ has to be almost the support of $\Pi_n$:  it is  guaranteed by Assumption {\rm \textbf{H4}}, which  in addition quantifies how $\Pi_n$ covers ${\cal Q}_n$. If {\rm \textbf{H2}} holds true with  $B_{d_{\nu*}}(q_0,\epsilon)$ instead of $B_{d_{\nu*},n}(q_0,\epsilon)$, then ${\cal Q}_n$ coincides with ${\cal Q}$ and Assumption {\rm \textbf{H4}} is no more needed.

\item Although  our semi-Markov framework differs from the Markov one, it is worth noticing that  Assumption  {\rm \textbf{H1}} is similar to the  one stated as Equation (4.1)
 in Ghosal and van Der Vaart \cite{GDV07}. In particular, for Markov chains,  this assumption is related to the transition probabilities of the Markov chain, whereas in our context,  {\rm \textbf{H1}} is concerned with the SMK density.



%
\end{itemize}

Note also that Assumption $\textbf{H1}$ could be replaced by the following:


\begin{itemize}
\item {\rm $\widetilde{\textbf{H1}}$}: 
There exists a  strictly positive constant $C$ and a strictly positive integer $k$ such that for any ${x}\in{E}$,
\begin{equation*}
{\frac{1}{k}\sum_{u=1}^{k}{P^{(u)}{(x,\cdot)}}}  \geq C.
\end{equation*}
\end{itemize}

\section{Proofs}

\subsection{Proof of Proposition \ref{prop:caracterisation}}\label{subsec:proof-prop:caracterisation}

In order to prove Proposition \ref{prop:caracterisation}, we prove that the right-hand side of Eq (\ref{eq:carac}) satisfies the two relevant
conditions. First,  for any $A \in \mathcal{E}$, any  $ \Gamma \in \mathbb{B}^+$, we have
\begin{eqnarray*}
\widetilde{\boldsymbol{\rho}}{Q(A,{\Gamma})}&:=&\int_{E\times{\mathbb{R}^{+}}}^{}{\widetilde{\boldsymbol{\rho}}(dy,{ds})}Q_y(A,\Gamma)\\
&=&\int_{E\times{E}\times{\mathbb{R}^{+}}}{\boldsymbol{\rho}(dx)Q_x (dy,{ds})}Q_y(A,{\Gamma})\\
&=&\int_{E}^{}{\boldsymbol{\rho}(dy)Q_y(A,{\Gamma})}\\
&=&\widetilde{\boldsymbol{\rho}}(A,{\Gamma}).
\end{eqnarray*}

Second,
\begin{eqnarray*}
\widetilde{\boldsymbol{\rho}}(E,{\mathbb{R}^{+}})=\int_{E}^{}{\boldsymbol{\rho}(dx)Q_x(E,{\mathbb{R}^{+}})}=1.
\end{eqnarray*}

\subsection{Proof of Proposition \ref{Prop:test}}\label{subsec:proof-Prop}
Our  proof is  constructive; indeed, we are going to construct a suitable testing procedure, namely $\psi_1({\cal H}_n)$, for the hypotheses testing problem given in  \eqref{eq_HypoTest}, i.e., 
\begin{eqnarray*}
H_{0}: q_{0} \;\text{against} \; H_{1}: q \in B_{d_{\eta*},n}(q_{1},\xi\epsilon), \; \text{{\rm with}} \;d_{\nu*}(q_{0},q_{1})\geq{\epsilon}, \; \mbox{{\rm and some } }    \xi \in (0,1).
\end{eqnarray*}

%
To control  exponentially both the type I and type II errors of  $\psi_1({\cal H}_n)$, we first  fix  some $x\in E$ for which we construct  the ``least favorable'' pair of RN derivatives of semi-Markov kernels  associated to the  following auxiliary testing problem 
\begin{align}
{\widetilde H}_{0, x}  :  q_{x;0}(\cdot,\cdot) \;   \mbox{{\rm against}} \;   {\widetilde H}_{1, x}  :\;  \left\{q_{x}(\cdot,\cdot) : h^2_\nu (Q_{x},Q_{x;1})\leq 1-\cos(\lambda \alpha_x)\right\},\label{pb-test-princ}  
\end{align}
where   $\lambda$ is any value in $]0,1/4[$ and $\alpha_x$ belongs to $]0,\pi/2[$ such that  
\begin{align}
h_{\nu}^2 ( Q_{x;0} ,Q_{x;1})= 1 -\cos (\alpha_x). \label{hellin-1-0}
\end{align}  

Based on this least favorable pair of $q_x$'s, we will then derive the construction of   $\psi_1({\cal H}_n)$ for the testing problem \eqref{eq_HypoTest}. 

For the sake of simplicity, let us denote  by  $q_{x}$ and $q_{x;j}$ for $j \in \BBn$ the probability density functions $q_{x}(\cdot,\cdot)$ and  $q_{x;j}(\cdot,\cdot)$, respectively.    

\subsubsection*{Least favorable pair of $q_x$'s for the testing problem  \eqref{pb-test-princ}}
For our purposes, we adapt the construction of Birg\'e \cite{B83a} for Markov chains to the semi-Markov framework. 
Whatever  is  $x$ in ${E}$, we attach to $x$ a particular probability density function $q_{x;2} \in  {\tilde H}_{1,x}$ defined by  
$$ q_{x;2}= \left(\frac{\sin((1-\lambda)\alpha_x)}{\sin(\alpha_x)}   \sqrt{ q_{x;1}} + \frac{\sin ( \lambda \alpha_x)}{\sin(\alpha_x)} \sqrt{q_{x;0}}\right)^2. $$ 
By construction, the following relations hold: 
\begin{eqnarray}
\lambda^2  h^2_\nu (Q_{x;0},Q_{x;1}) & \leq &  h^2_\nu (Q_{x;1},Q_{x;2})   \leq  h^2_\nu (Q_{x;0},Q_{x;1});\label{trick-helli-1}  \\
 (1-\lambda)^2 h^2_\nu (Q_{x;0},Q_{x;1}) &\leq & h^2_\nu (Q_{x;0},Q_{x;2});
 \label{trick-helli-2}\\
  h^2_\nu (Q_{x;1},Q_{x;2}) & =& 1 - \cos (\lambda \alpha_x); \label{trick-helli-3}\\
   h^2_\nu (Q_{x;0},Q_{x;2}) & =& 1 - \cos ((1-\lambda) \alpha_x). \nonumber 
\end{eqnarray}

\subsubsection*{Construction of the test procedure for  the testing problem \eqref{eq_HypoTest}}\label{subsubsec:construc-test}

\vspace{0.3cm}
Next, we set $\kappa=k+l$ and $N=[n/{\kappa}]$, where $l$ and $k$ are issued from Assumption {\bf H1} and $\big{[}\cdot\big{]}$ denotes
the integer part.  We consider $N$ i.i.d. random variables
$Y_{1},Y_{2},\ldots,Y_{N}$, which are generated independently from $\mathcal{H}_{n}$ according to the discrete uniform distribution $\mathcal{U}_{\{ 1,\ldots,k\}}$. 

We further define the test statistic
\begin{eqnarray*}
T(\mathcal{H}_{n})&=&\sum_{i=1}^{N}\log{\Phi_{J_{\tau_{i}-1}}(J_{\tau_{i}},X_{\tau_{i}})},
\end{eqnarray*}
where
$$\left\{ \begin{array}{rcl}
\Phi_{J_{\tau_{i}-1}}(J_{\tau_{i}},X_{\tau_{i}})&=&\sqrt{\frac{q_{J_{\tau_{i}-1};2}(J_{\tau_{i}},X_{\tau_{i}})}{q_{J_{\tau_{i}-1};0}(J_{\tau_{i}},X_{\tau_{i}})}},\\
\tau_{i } & =& \kappa(i-1)+l+Y_{i}.\end{array}\right.$$

Our test procedure for the hypotheses problem \eqref{eq_HypoTest}  is then defined as follows
\begin{equation} \label{test-stat}
\psi_{1} (\mathcal{H}_{n})=\BBone_{\{T(\mathcal{H}_{n})>0 \} }.
\end{equation}
%



\item \textsc{Test simple hypothesis vs simple hypothesis}

Let us focus on the general SMPs and consider the
following statistical test:
\begin{eqnarray*}
H_{0}: q_{0} \;\text{against} \; H_{1}: q_{1} \; \text{{\rm with}} \;d_{\nu*}(q_{0},q_{1})\geq{\epsilon}.
\end{eqnarray*}

To construct the testing procedure, the test statistic defined in (\ref{test-stat}), should
be modified as follows:
\begin{eqnarray*}
T(\mathcal{H}_{n})&=&\sum_{i=1}^{N}\log{\Phi_{J_{\tau_{i}-1}}(J_{\tau_{i}},X_{\tau_{i}})},
\end{eqnarray*}
where
$$\left\{ \begin{array}{rcl}
\Phi_{J_{\tau_{i}-1}}(J_{\tau_{i}},X_{\tau_{i}})&=&\sqrt{\frac{q_{J_{\tau_{i}-1};1}(J_{\tau_{i}},X_{\tau_{i}})}{q_{J_{\tau_{i}-1};0}(J_{\tau_{i}},X_{\tau_{i}})}},\\
\tau_{i } & =& \kappa(i-1)+1+Y_{i}\\
\kappa & =& k+1\\
Y_i&\stackrel{iid}{\sim}& {\cal U}_{\{1,\ldots,k \}}.\end{array}\right.$$

In this case, Hypothesis {\bf H1} reduces to $\textbf{H1}^{\sharp}:$

\begin{itemize}
\item {\rm $\textbf{H1}^{\sharp}$}: 
There exist a measure $\nu^{*}$ on
$\mathcal{E}$ and a positive integer $k$ such that for any ${x}\in{E}$,

\begin{equation*}
{\frac{1}{k}\sum_{u=1}^{k}{P^{(u)}{(x,\cdot)}}}  \geq \nu^{*}(\cdot).
\end{equation*}

\end{itemize}

Then following the steps of the proof of the Proposition \ref{Prop:test}
and replacing the Assumption \rm $\textbf{H1}$ by {\rm $\textbf{H1}^{\sharp}$}
lead us to the desired result. It is worth mentioning that in this case
the inequalities
(\ref{trick-helli-1}), (\ref{trick-helli-2}), (\ref{trick-helli-3}) and
Lemma \ref{Lem:test} are not used. 

Note also that in Proposition \ref{Prop:test}, the upper-bound of both errors is the same, equal to  $\exp(-Kn\epsilon^2).$

\end{itemize}

\subsubsection*{Type I error probability}\label{subsubsec:typeI}
By means of the Markov property we obtain that
\begin{eqnarray}
\mathbb{E}_{0} (\psi_1({\cal H}_n)) 
&\leq&{\mathbb{E}_{0} 
\bigl{(}\prod_{i=1}^{N-1}{\Phi_{J_{\tau_{i}-1}}(J_{\tau_{i}},X_{\tau_{i}})}\Phi_{J_{\tau_{N}-1}}(J_{\tau_{N}},X_{\tau_{N}})\bigr{)}} \nonumber \\
&=&\mathbb{E}_{0}
\bigl{(}\prod_{i=1}^{N-1}{\Phi_{J_{\tau_{i}-1}}(J_{\tau_{i}},X_{\tau_{i}})}\mathbb{E}_{0}(\Phi_{J_{\tau_{N}-1}}(J_{\tau_{N}},X_{\tau_{N}})|{\mathcal{H}_{\kappa(N-1)}})\bigr{)} \nonumber \\
&=&\mathbb{E}_{0}
\bigl{(}\prod_{i=1}^{N-1}{\Phi_{J_{\tau_{i}-1}}(J_{\tau_{i}},X_{\tau_{i}})}\mathbb{E}_{0}({\Phi_{J_{\tau_{N}-1}}(J_{\tau_{N}},X_{\tau_{N}})}|J_{\kappa(N-1)})\bigr{)},\label{1espece-global}
\end{eqnarray}
where
$\mathcal{H}_{\kappa(N-1)}=(J_{0},\ldots,J_{\kappa(N-1)}, X_0,\ldots,X_{\kappa(N-1)},).$

\medskip
\begin{itemize}
\item \textbf{Step 1} 

Set $T_1 := \mathbb{E}_{0}({\Phi_{J_{\tau_{N}-1}}(J_{\tau_{N}},X_{\tau_{N}})}|J_{\kappa(N-1)})$.  Since  $\tau_{i}\sim{U_{\{\kappa(i-1)+l+1,\ldots,\kappa{i}\}}}$, we obtain
\begin{eqnarray*}
T_{1}&=&\frac{1}{k}\sum_{u=1}^{k}{\mathbb{E}_{0}\big{[}{\Phi_{J_{\kappa(N-1)+l+u-1}}{(J_{\kappa(N-1)+l+u},X_{\kappa(N-1)+l+u})} \vert J_{\kappa(N-1)}}\big{]}}.
\end{eqnarray*}
Next set $\Gamma_{u}:=\mathbb{E}_{0}\big{[}{\Phi_{J_{\kappa(N-1)+l+u-1}}{(J_{\kappa(N-1)+l+u},X_{\kappa(N-1)+l+u})} \vert J_{\kappa(N-1)}}\big{]}$ and rewrite $\Gamma_u$ as follows, 
\begin{eqnarray*}
\Gamma_{u}&=&\int_{E}^{}\int_{E}^{}\int_{\mathbb{R}^{+}}^{}{\Phi_{x}(y,t)}P^{(l+u-1)}_{0}(J_{\kappa(N-1)},dx)
q_{x;0}(y,t)d{\nu(y,t)}\\
&=&\int_{E}^{}P^{(l+u-1)}_{0}(J_{\kappa(N-1)},dx)\int_{E}^{}\int_{\mathbb{R}^{+}}^{}{\Phi_{x}(y,t)}
q_{x;0}(y,t)d{\nu(y,t)}\\
&=&\int_{E}^{}P^{(l+u-1)}_{0}(J_{\kappa(N-1)},dx){\Big{(}1-h_{\nu}^{2}(Q_{x;0},Q_{x;2})\Big{)}},
\end{eqnarray*}
where the last equality is due to $$\int_{E} \int_{\mathbb{R}^{+}}  \sqrt{q_{x;2}q_{x;0} }d\nu  =1 - h^2_{\nu} (Q_{x;0},Q_{x;2}).$$

Assumption \textbf{H1} and Eq. (\ref{trick-helli-2}) lead us to the following upper bound of $T_{1}$:
\begin{eqnarray*}
T_{1}&=&1-\frac{1}{k}\sum_{u=1}^{k}\int_{E}^{}P^{(l+u-1)}_{0}(J_{\kappa(N-1)},dx)h_{\nu}^{2}
(Q_{x;0},Q_{x;2})\\
&\leq &{1-\int_{E}^{}h_{\nu}^{2}(Q_{x;0},Q_{x;2})d\nu^{*}(x)}\\
& \leq & 1- (1-\lambda)^2\int_{E}^{} h^2_{\nu}(Q_{x;0},Q_{x;1})  d\nu^{*}(x)\\
&=&1-(1-\lambda)^2 d^{2}_{\nu^{*}}(q_{0},q_{1})\\
&\leq&{e^{- (1-\lambda)^2d^{2}_{\nu^{*}}(q_{0},q_{1})}}\\
&\leq&{e^{-(1-\lambda)^{2}\epsilon^{2}}}.
\end{eqnarray*}
This latter inequality provides  a first upper bound of $\mathbb{E}_{0} ( \psi_{1}(\mathcal{H}_{n}))$ via the relation \eqref{1espece-global}.

\item Then, by setting $T_i   := \mathbb{E}_{0}({\Phi_{J_{\tau_{N-i+1}-1}}(J_{\tau_{N-i+1}},X_{\tau_{N-i+1}})}|J_{\kappa(N-i)})$ for $i=2,\ldots,N$, 
and by repeating \textbf{Step 1} for the successive $T_i$, we finally obtain
\begin{eqnarray*}
\mathbb{E}_{0}
\bigl{(}\psi_{1}(\mathcal{H}_{n})\bigr{)}&\leq&{e^{-\frac{n}{\kappa}(1-\lambda)^2\epsilon^{2}}}
=e^{-Kn\epsilon^{2}}, \quad \mbox{ {\rm with} } K=\frac{(1-\lambda)^2}{\kappa}.
\end{eqnarray*}
\end{itemize}






\subsubsection*{Type II error probability}\label{subsubsec:typeII}

To bound from above the  type II error probability, we need an additional result stated as Lemma \ref{Lem:test}. This lemma provides upper bounds  for a  quantity which is similar to the $T_1$-term appearing in the first type error probability. The main difference here is that this quantity should be bounded from above uniformly over $q $ in $ B_{d_{\eta^*},n} (q_1, \xi\epsilon)$.

This requires the definition of the subset $G_q$ of 
$E$ by
\begin{eqnarray*}
G_q:
=\{x\in{E}:h_{\nu}(Q_{x},Q_{x;1})\leq{\lambda{h_{\nu}(Q_{x;0},Q_{x;1})}}\},
\end{eqnarray*}
and the notation of its complementary into $E$ by $G_q^{\complement}$.

\begin{Lemma}\label{Lem:test}
For any $\lambda \in ]0,1/4[$, there exists $\iota \in [0,\frac{3}{4}[$, such that for all $q \in  B_{d_{\eta^*},n} (q_1, \xi\epsilon)$, 

\begin{itemize}
\item if $x\in{G_q}$, then
\begin{eqnarray}\label{eq:U}
\mathbb{E}_q[\Phi^{-1}_{J_0}(J_1,X_1)\vert J_0=x]\leq{1-h_{\nu}^{2}(Q_{x;{0}},Q_{x;2})}\leq{1-(1-\lambda)^{2}h_{\nu}^{2}(Q_{x;0},Q_{x;1})};
\end{eqnarray}

\item if $x\in{G_q^{\complement}}$, then
\begin{eqnarray}
\mathbb{E}_{q}[\Phi^{-1}_{J_0}(J_1,X_1)\vert J_0=x] & < & 1   +    8 \frac{1-\lambda}{\lambda} h_\nu^2 (Q_{x},Q_{x;1}) \nonumber \\ 
&- & (1-\frac{2 \lambda}{1-\lambda}) [   1 - \iota]    h^2_\nu  ( Q_{x;0},Q_{x;1}). 
\label{eq:complU}
\end{eqnarray} 
\end{itemize}
\end{Lemma}
The proof of Lemma \ref{Lem:test} is postponed to Section  \ref{subsubsec:Lemma-Test}.  
%

Consider $\Phi^{-1}$ equal to one over $\Phi$,  that is  $\displaystyle{\Phi^{-1} =\sqrt{\frac{q_{0}}{q_{2}}}}$.  Similarly to the calculations of the type I error probability, we obtain that for any 
$ q \in B_{d_{\eta*},n} (q_{1}, \xi \epsilon)$,
\begin{eqnarray*}
\mathbb{E}_q
\bigl{(}1-\psi_{1}(\mathcal{H}_{n})\bigr{)}  
&\leq &\mathbb{E}_q
\big{(}\prod_{i=1}^{N-1}{\Phi^{-1}_{J_{\tau_{i}-1}}(J_{\tau_{i}},X_{\tau_{i}})}\mathbb{E}_q({\Phi^{-1}_{J_{\tau_{N}-1}}(J_{\tau_{N}},X_{\tau_{N}})}|J_{\kappa(N-1)})\big{)}.
\end{eqnarray*}


Similarly to $T_1$,  we further define $W_{1}$  by
\begin{eqnarray*}
W_1&:=&\mathbb{E}_q({\Phi^{-1}_{J_{\tau_{N}-1}}(J_{\tau_{N}},X_{\tau_{N}})}|J_{\kappa(N-1)})\\
& =& \frac{1}{k}\sum_{u=1}^{k}{\mathbb{E}_q\big{[}{\Phi^{-1}_{J_{\kappa(N-1)+l+u-1}}{(J_{\kappa(N-1)+l+u},X_{\kappa(N-1)+l+u})}|J_{\kappa(N-1)}}\big{]}}.
\end{eqnarray*}

\begin{itemize}
\item \textbf{Step 2} 
Taking into account the partition of $E$ into $G_q$ and $G_q^{\complement}$, we obtain
\begin{eqnarray*}
W_{1}&=&\frac{1}{k}\sum_{u=1}^{k} \int_{\mathbb{R}^{+}}^{}\int_{E}^{}\int_{E}^{}{\Phi^{-1}_{x}(y,t)}P_q^{(l+u-1)}(J_{\kappa(N-1)},dx)q_x(y,t)d{\nu(y,t)}\\
&=&\frac{1}{k}\sum_{u=1}^{k} \int_{E}^{}P_q^{(l+u-1)}(J_{\kappa(N-1)},dx)\mathbb{E}_q[\Phi^{-1}_{J_0}(J_1,X_1)\vert J_0=x]\\
&=&\frac{1}{k}\sum_{u=1}^{k} \int_{G_q}^{}P_q^{(l+u-1)}(J_{\kappa(N-1)},dx)\mathbb{E}_q[\Phi^{-1}_{J_0}(J_1,X_1)\vert J_0=x]\\
&\quad &+\frac{1}{k}\sum_{u=1}^{k}\int_{G_q^{\complement}}^{}P_q^{(l+u-1)}(J_{\kappa(N-1)},dx)\mathbb{E}_q[\Phi^{-1}_{J_0}(J_1,X_1)\vert J_0=x]. 
\end{eqnarray*}

Combining with $(1-\lambda)^2 > \displaystyle{\frac{1 -3 \lambda}{1-\lambda}}$ , Assumption \textbf{H1} and Lemma \ref{Lem:test} lead to,
\begin{eqnarray*}
W_1&\leq& 1 - \frac{1 -3 \lambda}{1-\lambda}  [   1 - \iota]  \frac{1}{k} \sum_{u=1}^{k}  \int_{E} P_q^{(l+u-1)}(J_{\kappa(N-1)},dx)h_{\nu}^{2}(Q_{x;0},Q_{x;1})  \nonumber\\
&+ & 8 \frac{1-\lambda}{\lambda} \frac{1}{k} \sum_{u=1}^{k}\int_{G_q^{\complement}} P_q^{(l+u-1)}(J_{\kappa(N-1)},dx)    h_{\nu}^{2}(Q_{x},Q_{x;1})  \nonumber \\
& \leq & 1 - \frac{1 -3 \lambda}{1-\lambda}  [   1 - \iota]  \int_{E}   h_{\nu}^{2}(Q_{x;0},Q_{x;1}) d\nu^*(x) +8 \frac{1-\lambda}{\lambda}
\int_{E}   h_{\nu}^{2}(Q_{x},Q_{x;1}) d\eta^*(x)\nonumber \\
& =& 1 - \frac{1 -3 \lambda}{1-\lambda}  [   1 - \iota]    d^2_{\nu^*} (q_{0},q_{1}) +
8 \frac{1-\lambda}{\lambda} d^2_{\eta^*} (q,q_{1}) \nonumber \\
&\leq & \exp \left( - \left\{ \frac{1 -3 \lambda}{1-\lambda}  [   1 - \iota] - 8 \frac{1-\lambda}{\lambda}\xi^2
\right\}\epsilon^2\right)= \exp \left( - K(\lambda) 
\epsilon^2\right),
\end{eqnarray*}
where $K(\lambda)$  is positive since there exists $\xi >0$ such that $\displaystyle{\frac{1 -3 \lambda}{1-\lambda}  [   1 - \iota] }> 8 \displaystyle{{\frac{(1-\lambda)}{\lambda}\xi^2}}$.

\item To complete the proof, we consider $W_i   := \mathbb{E}_{q}({\Phi^{-1}_{J_{\tau_{N-i+1}-1}}(J_{\tau_{N-i+1}},X_{\tau_{N-i+1}})}|J_{\kappa(N-i)})$ for $i=2,\ldots,N$. We then repeat \textbf{Step 2} for the successive $W_i$, and finally deduce that for any $q \in B_{d_{\eta*},n} (q_{1}, \xi \epsilon)$,
\begin{eqnarray*}
\mathbb{E}_q^{(n)}\bigl{(}1-\psi_{1}(\mathcal{H}_{n})\bigr{)} & \leq & \exp\Big(-n \tilde{K}(\lambda) \epsilon^{2}\Big),
\end{eqnarray*}
with $\tilde{K}(\lambda)=K(\lambda)/\kappa$.
%
%
\end{itemize}

\subsection{Proof of Lemma \ref{Lem:test} }\label{subsubsec:Lemma-Test}
We define the Hellinger affinity between two distributions $P_1$ and $P_2$, absolutely continuous w.r.t. $\nu$ ,with derivatives $p_1$ and $p_2$ respectively, by 
\begin{eqnarray*}
\varrho_{\nu}(P_1,P_2)&:=&\int_{\mathbb{R}^{+}}\int_{E}\sqrt{p_1 p_2} d{\nu} 
=1-h_{\nu}^{2}(P_1,P_2).
\end{eqnarray*}

\medskip
In the sequel, let $q$ be an arbitrary element of $B_{d_{\eta^*},n} (q_{1}, \xi \epsilon)$.

When $x$ belongs to $G_q$,  the proof of \eqref{eq:U} results  directly from Theorem 2 in Birg\'e \cite{B2013}.

\bigskip
When   $x$ belongs to $G_q^{\complement}$, i.e., $x \in E$ such that  
$h_{\nu}(Q_{x},Q_{x;1})> \lambda{h_{\nu}(Q_{x;0},Q_{x;1})}$, let us  prove the statement (\ref{eq:complU}). 
 
We follow the lines of Birg\'e \cite{B83a} and consider a real number $A$ such that $A\geq \displaystyle{\frac{2}{1 -\lambda}}$. We then decompose the 
term $\mathbb{E}_q[\Phi^{-1}_{J_0}(J_1,X_1)\vert J_0=x]$ into four terms: 
\begin{eqnarray*}
\mathbb{E}_q[\Phi^{-1}_{J_0}(J_1,X_1)\vert J_0=x] & \leq & \mathbb{E}_{q_1}[\Phi^{-1}_{J_0}(J_1,X_1)\vert J_0=x] 
+ \sum_{i=1}^3 \int_{{\cal A}_{x;i}} (\Phi^{-1}_{x} -1) (q_{x}-q_{x;1}) d \nu\\
&  :=&  T_0 + \sum_{i=1}^3 T_i,
\end{eqnarray*}
where 
\begin{align*}
{\cal A}_{x;1} =& \left\{(y,t)\in{E\times\mathbb{R}^{+}}: \sqrt{\frac{q_{x}(y,t)}{q_{x;1}(y,t)}} > A-1, \text{ } \Phi^{-1}_{x}(y,t) >1    \right\}\\
{\cal A}_{x;2} = & \left\{(y,t)\in{E\times\mathbb{R}^{+}}: 1 \leq  \sqrt{\frac{q_{x}(y,t)}{q_{x;1}(y,t)}} \leq  A-1, \text{ } \Phi^{-1}_{x}(y,t) >1    \right\}\\
{\cal A}_{x;3} =&  \left\{(y,t)\in{E\times\mathbb{R}^{+}}: \sqrt{\frac{q_{x}(y,t)}{q_{x;1}(y,t)}} < 1, \text{ } \Phi^{-1}_{x}(y,t) <1  \right\}.
\end{align*}

For the sake of simplicity, set $r_x=\displaystyle{\frac{q_{x}}{q_{x;1}}}$ and start with $T_0$. Due to the definition  of $\Phi^{-1}_{x}(\cdot,\cdot)$, to Equation  \eqref{hellin-1-0} and to the concavity of the function $\displaystyle{y \rightarrow \frac{\sin(\alpha_x)y}{\sin(\alpha_x\lambda ) y + \sin(\alpha_x(1-\lambda ))}}$ , we deduce that
 \begin{eqnarray}
 T_{0}  & \leq & \frac{\sin(\alpha_x) \rho_\nu (Q_{x;0},Q_{x;1})}{\sin(\alpha_x\lambda ) \rho_\nu (Q_{x;0},Q_{x;1}) + \sin(\alpha_x(1-\lambda ))}\nonumber \\
 & =& \frac{\sin(\alpha_x)  \cos (\alpha_x)}{\sin(\alpha_x\lambda )  \cos (\alpha_x) + \sin(\alpha_x(1-\lambda ))}\nonumber \\
 & =& \frac{ \cos (\alpha_x) }{ \cos (\alpha_x\lambda) } \nonumber \\
 & \leq & 1-\biggl{(}1-\frac{2 \lambda}{1-\lambda}\biggr{)} h^2_\nu ( Q_{x;0},Q_{x;1}) ,\label{Lem:T_0}
 \end{eqnarray}
 where the last inequality results from both the convexity of the $tan$ function on $]0,\pi/2[$  and $\lambda < 1/4$. 
 
 \medskip
 Let us now turn to $T_1$. First note that 
 \begin{eqnarray}
 \Phi^{-1}_{x} 
& =& \frac{\sin(\alpha_x)  \sqrt{\frac{q_{x;0}}{q_{x;1}}}}{\sin(\alpha_x\lambda ) 
 \sqrt{\frac{q_{x;0}}{q_{x;1}}} + \sin(\alpha_x(1-\lambda ))} \nonumber  \\
 & \leq & \frac{\sin(\alpha_x)  }{\sin(\alpha_x\lambda ) 
  } < \frac{1}{\lambda},\label{phi-1}
 \end{eqnarray} 
 where (\ref{phi-1}) results from the following inequality 
 \begin{eqnarray}
\forall\;  \alpha \in ]0, \pi/2[,\quad \forall \lambda \in ] 0, 1[, \quad  \frac{\sin(\lambda \alpha)}{\lambda \sin (\alpha)}>1. \label{use-ine-sin}
\end{eqnarray}

On ${\cal A}_{x;1}$, since $\displaystyle{r_x -1 <\frac{A}{A-2} (\sqrt{r_x} -1)^2}$, then from (\ref{phi-1}) we obtain,

 \begin{eqnarray}
 T_1  
 & \leq & \frac{A}{A-2} \frac{1-\lambda}{\lambda} \int_{{\cal A}_{x;1}}  (\sqrt{q_x }-\sqrt{ q_{x;1}} )^2 d\nu \nonumber \\
  & \leq &      \frac{A}{A-2} \frac{1-\lambda}{\lambda}   2 h^2_{\nu} (Q_x,Q_{x;1}) - \frac{A}{A-2} \frac{1-\lambda}{\lambda} \int_{{\cal A}_x }   (\sqrt{q_x }-\sqrt{ q_{x;1}} )^2 d\nu,   \label{Lem:T_1}
  \end{eqnarray}
where ${\cal A}_x$ is a subset of ${{\cal A}_{x;1}^{\complement}}$.

Second we study the last two terms $T_2$ and $T_3$. 
On ${\cal A}_{x;2}$ and ${\cal A}_{x;3}$, we first apply the Cauchy-Schwarz inequality,  i.e., $\forall i \in \{2,3\}$,
\begin{eqnarray*}
\left(\int_{{\cal A}_{x;i}} (\Phi^{-1}_{x} -1) (r_x-1)  q_{x;1} d \nu  \right)^2
\leq  \int_{{\cal A}_{x;i}} (\Phi^{-1}_{x} -1)^2   q_{x;1} d \nu  \int_{{\cal A}_{x;i}}(r_x-1)^2  q_{x;1} d \nu.
\end{eqnarray*}
Second we note that 
\begin{eqnarray}
\int_{{\cal A}_{x;i}} (\Phi^{-1}_{x}(\cdot,\cdot) -1)^2   q_{x;1} d \nu
& =& \int_{{\cal A}_{x;i}} (\sqrt{q_{x;0}} - \sqrt{q_{x;2}})^2 
\frac{q_{x;1}}{q_{x;2}} d \nu  \nonumber \\
& \leq & \beta \int_{{\cal A}_{x;i}} (\sqrt{q_{x;0}} - \sqrt{q_{x;2}})^2  d \nu \label{term=T21T31},
\end{eqnarray}
where $\beta$, the upper bound of $\displaystyle{\frac{q_{x;1}}{q_{x;2}}}$,  
is given by 
$\beta =\left\{ \begin{array}{ccll} 1 & \mbox{{\rm on }} & {\cal A}_{x;2} & \mbox{{\rm since} } \displaystyle{\frac{q_{x;1}}{q_{x;0}}} <1,\\
\displaystyle{\frac{1}{(1-\lambda)^2}} & \mbox{{\rm on }} & {\cal A}_{x;3} & \mbox{{\rm due to \eqref{use-ine-sin}}}.
\end{array} \right. $ 
We further note that 
$$   \int_{{\cal A}_{x;i}}(r_x-1)^2  q_{x;1}(\cdot,\cdot) d \nu  \leq \left\{ \begin{array}{c}  
A^2 \int_{{\cal A}_{x;2}} (\sqrt{q_x} - \sqrt{q_{x;1}})^2  d \nu   \\
2^2 \int_{{\cal A}_{x;3}} (\sqrt{q_x} - \sqrt{q_{x;1}})^2  d \nu
\end{array}
.
\right.     $$

The latter combined with (\ref{term=T21T31}) and since $A > 2/(1-\lambda)$, entails  
\begin{eqnarray}
T_2 + T_3  \leq &  A \left( \int_{{\cal A}_{x;2}} (\sqrt{q_x} - \sqrt{q_{x;1}})^2  d \nu   \int_{{\cal A}_{x;2}} (\sqrt{q_{x;0}} - \sqrt{q_{x;2}})^2  d \nu\right)^{1/2} \nonumber \\ 
 &+ \displaystyle{\frac{2}{1-\lambda}} \left( \int_{{\cal A}_{x;3}} (\sqrt{q_x} - \sqrt{q_{x;1}})^2  d \nu  \int_{{\cal A}_{x;3}} (\sqrt{q_{x;0}} - \sqrt{q_{x;2}})^2  d \nu \right)^{1/2} \nonumber \\
\leq &   A   \left( \int_{{\cal A}_{x}} (\sqrt{q_x} - \sqrt{q_{x;1}})^2  d \nu  
\int_{{\cal A}_{x;2} \cup  {\cal A}_{x;3}} (\sqrt{q_{x;0}} - \sqrt{q_{x;2}})^2  d \nu\right)^{1/2}. \label{T2T3}
\end{eqnarray}

From (\ref{Lem:T_0}), (\ref{Lem:T_1}) and (\ref{T2T3}), it follows that
\begin{eqnarray}
 \mathbb{E}[\Phi^{-1}_{J_0}(J_1,X_1)\vert J_0=x]  & \leq &  1-\biggl{(}1-\frac{2 \lambda}{1-\lambda}\biggr{)} h^2_\nu ( Q_{x;0},Q_{x;1})   + 2\frac{A}{A-2} \frac{1-\lambda}{\lambda} h_\nu^2 (Q_{x},Q_{x;1}) \nonumber \\
& &
  -  \frac{A}{A-2} \frac{1-\lambda}{\lambda} \int_{{\cal A}_{x} }   (\sqrt{q_x }-\sqrt{ q_{x;1}} )^2 d\nu \nonumber \\ 
 & &  + A   \left( \int_{{\cal A}_{x} } (\sqrt{q_x} - \sqrt{q_{x;1}})^2  d \nu  
\int_{{\cal A}_{x;2} \cup  {\cal A}_{x;3}} (\sqrt{q_{x;0}} - \sqrt{q_{x;2}})^2  d \nu\right)^{1/2}. \nonumber 
\end{eqnarray}
At a next step we consider the following function of $z_x$  
\begin{equation*}
z_x\rightarrow   -  \frac{A}{A-2} \frac{1-\lambda}{\lambda} z_x 
+ z_x^{1/2}A \left(\int_{{\cal A}_{x;2} \cup  {\cal A}_{x;3}} (\sqrt{q_{x;0}} - \sqrt{q_{x;2}})^2  d \nu\right)^{1/2},
\end{equation*}
whose maximum is reached at  
$$z_{x;max}= \frac{1}{4}(A-2)^2 \biggl{(}\frac{\lambda}{1-\lambda}\biggr{)}^2 \int_{{\cal A}_{x;2} \cup  {\cal A}_{x;3}} (\sqrt{q_{x;0}} - \sqrt{q_{x;2}})^2  d \nu. $$

Hence we obtain a new upper bound of $\mathbb{E}_{1}[\Phi^{-1}_{J_0}(J_1,X_1)\vert J_0=x]$, 
that is
\begin{eqnarray}
 \mathbb{E}_{q}[\Phi^{-1}_{J_0}(J_1,X_1)\vert J_0=x]  
 & \leq &  1-\biggl{(}1-\frac{2 \lambda}{1-\lambda}\biggr{)} h^2_\nu ( Q_{x;0},Q_{x;1})   + 2\frac{A}{A-2} \frac{1-\lambda}{\lambda} h_\nu^2 (Q_{x},Q_{x;1}) \nonumber \\
& &
 + \frac{A(A-2)}{2} \frac{\lambda}{1-\lambda}  h^2_\nu (Q_{x;0},Q_{x;2})\nonumber \\
 & \leq & 
 1   +    2\frac{A}{A-2} \frac{1-\lambda}{\lambda} h_\nu^2 (Q_{x},Q_{x;1})  -(1-\frac{2 \lambda}{1-\lambda}) h^2_\nu ( Q_{x;0},Q_{x;1}) \nonumber \\ 
&& + A(A-2) \frac{\lambda}{1-\lambda}  h^2_\nu  ( Q_{x;0},Q_{x;1})  \sin^2 \Big{(}(1-\lambda)\frac{\pi}{4}\Big{)}, \nonumber \\
& \leq & 
   1   +     2\frac{A}{A-2} \frac{1-\lambda}{\lambda} h_\nu^2 (Q_{x},Q_{x;1}) \nonumber \\ 
& &  - \biggl{(}1-\frac{2 \lambda}{1-\lambda}\biggr{)} \Big{[}   1 - \frac{A(A-2)\lambda}{(1-3\lambda)} \sin^2 \Big{(}(1-\lambda)\frac{\pi}{4}\Big{)}   \Big{]}    h^2_\nu  ( Q_{x;0},Q_{x;1}), \nonumber 
  \end{eqnarray}
where the penultimate inequality results from the increase of the function  $x\in ]0,\pi/2[ \rightarrow \displaystyle{\frac{\sin(\lambda x/2)}{\lambda  \sin(x/2)}}$
 for any $\lambda \in ]0,1]$. 
 
 \medskip
Finally, by setting $A= 8/3$ that satisfies $A \geq 2/(1-\lambda)$ and using both inequalities  $\displaystyle{\sin^2 \Big{(}(1-\lambda)\frac{\pi}{4}\Big{)}}< (1-\lambda)^2 \Big{(}\frac{\pi}{4}\Big{)}^2$ $\forall \lambda \in ]0,1/4[$  and  
 $\displaystyle{\frac{\lambda (1- \lambda)^2}{1-3\lambda}  <  \frac{9}{16}}$  $\forall \lambda \in ]0,1/4[$,
Lemma \ref{Lem:test} is proved with  $\iota=\frac{\pi^2}{16} <3/4$. \hfill $\Box$

\subsection{Proof of Corollary \ref{cor:test-global}}

The proof of Corollary \ref{cor:test-global}  is similar to the proof of Lemma $9$ in \cite{GDV07}. However, we sketch it in order to define the statistical test procedure $\psi({\cal H}_n)$. 
First, consider the partition: 
%
\begin{eqnarray*}
\{q \in {\cal Q}_n: d_{\nu*}( q_0,q) > \epsilon_{n}M \}& =& \bigcup_{j\geq{1}}\Big{\{}{q} \in{{{\mathcal{Q}}_n}}: j\epsilon_{n}M<d_{\nu*}({q}_0,{q})\leq{(j+1)\epsilon_{n}M}\Big{\}}\\
& =:& \bigcup_{j\geq{1}}{H}_{j}.
\end{eqnarray*}

For $\xi \in ]0,1[$, and any $j \geq 1$, we consider ${\widetilde{H}}_{j}$, a $j\epsilon_{n}\xi{M}$-net on ${H}_{j}$ for the distance $d_{\eta*}$  satisfying
three conditions:
\begin{itemize}
\item $\forall q \in {\widetilde{H}}_{j}$,  $d_{\nu*}(q_0,q)\geq  j\epsilon_{n}M$;
\item $\forall {q} \in{H}_{j}$, $\exists q_j \in {\widetilde{H}}_{j}$ such that  $d_{\eta*}(q ,q_j )\leq j\epsilon_{n}\xi{M}$;
\item  $\log{N}\bigl{(}{\epsilon_n M\xi},{\widetilde{H}}_{j},d_{\eta*}\bigr{)}\leq{n\epsilon_n^{2}}$ (due to \textbf{H2}).
\end{itemize}
 

For  $j \geq 1$ and any $q_{j,i} \in \widetilde{H}_{j}$, we then apply Proposition \ref{Prop:test} with $\epsilon=jM\epsilon_{n}$ and $q_1=q_{j,i}$; this  implies the existence of a statistical procedure  $\psi_{j,i}({\cal H}_n)$ 
that satisfies \eqref{test:1-2-error}. 
%

We then define our test procedure 
\begin{eqnarray}
\psi(\mathcal{H}_{n}):=\displaystyle{\max_{j\geq{1}}\max_{q_{j,i} \in\widetilde{H}_{j}}{\psi_{j,i}({\cal H}_n)}}. \label{test}\end{eqnarray}

We further combine Assumption  \textbf{H2} and Proposition \ref{Prop:test} to obtain for $M$ large enough
\begin{eqnarray*}
\mathbb{E}_{0}^{(n)}[\psi(\mathcal{H}_{n})]&\leq & \sum_{j=1}^{\infty}\sum_{q_{j,i}
\in\widetilde{H}_{j}}^{}
\mathbb{E}_{0}^{(n)}[\psi_{j,i}(\mathcal{H}_n)]\\
& \leq &  e^{  n \epsilon_n^2}  \frac{e^{-K n\epsilon^{2}_{n}M^{2}}}{1-e^{-K n\epsilon^{2}_{n}M^{2}}}
 \leq  e^{-Kn\epsilon^{2}_{n}M/2^{2}},
\end{eqnarray*}
and 
\begin{eqnarray*}
\sup_{q \in{\bigcup_{j\geq 1}{H}_{j}}}
{\mathbb{E}_{q}^{(n)}[1-\psi(\mathcal{H}_{n})]\leq{\sup_{j>1}e^{-{\tilde K} nj^{2}\epsilon^{2}_{n}M^{2}}}}\leq{
e^{-{\tilde K} n\epsilon^{2}_{n}M^{2}}}.
\end{eqnarray*}

\subsection{Proof of Theorem \ref{THM}} \label{sec:preuve-THM}
Let $M$ be a positive constant. We first decompose the right-hand side of \eqref{eq:resu} in two parts 
\begin{eqnarray}
\Pi_{n}^{\mathcal{H}_{n}}\bigl{(}B^{\complement}_{d_{\nu^*}}(q_{0},\epsilon_{n}M)\bigr{)} \!\!\!\!\! & =&  \!\!\!\!\! \Pi_{n}^{\mathcal{H}_{n}}\bigl{(}B^{\complement}_{d_{\nu^*}}(q_{0},\epsilon_{n}M)\cap \mathcal{Q}_{n}\bigr{)}
+\Pi_{n}^{\mathcal{H}_{n}}\bigl{(}B^{\complement}_{d_{\nu^*}}(q_{0},\epsilon_{n}M)\cap  \mathcal{Q}_{n}^{\complement}\bigr{)} \nonumber \\
\!\!\!\!\! &\!\!\!\!\! =:\!\!\!\!\! &  A_1 + A_2.
\label{proof:decomp}
\end{eqnarray}

In the sequel,  each term in the right-hand
side of (\ref{proof:decomp}) is separately bounded from above:  for $A_1$, we apply Corollory \ref{cor:test-global}, whereas to upper bound $A_2$ we use \textbf{H3} and \textbf{H4}.

%
%

First, let us focus on $A_1$. Recall that $\mathcal{L}_{n}(q)$, the likelihood function of the sample path $\mathcal{H}_{n}$ evaluated at $q \in {\cal Q}$,  is given by
\begin{eqnarray*}
\mathcal{L}_{n}(q) &=&\widetilde{\boldsymbol{\rho}}(J_{0},S_{0}){\prod_{l=1}^{n}q_{J_{l-1}}(J_{l},X_{l})}.
\end{eqnarray*}
 
Then, $A_1$ could be written as follows:
\begin{eqnarray*}
A_1
&=&\frac{\int_{B_{\nu*}^{\complement}(q_{0},\epsilon_{n}M)\cap\mathcal{Q}_n}^{}\mathcal{L}_{n}(q) d{\Pi_{n}}(q)}{\int_{\mathcal{Q}}^{}\mathcal{L}_{n}(q)d{\Pi_{n}}(q)} \\
&=&\frac{\int_{B_{\nu*}^{\complement}(q_{0},\epsilon_{n}M)\cap\mathcal{Q}_n} \frac{\mathcal{L}_{n}(q)}{\mathcal{L}_{n}(q_0)}d{\Pi_{n}}(q)}
{\int_{\mathcal{Q}}^{}\frac{\mathcal{L}_{n}(q)}{\mathcal{L}_{n}(q_0)}d{\Pi_{n}}(q)} \\ 
&:=&\frac{N_{n}}{D_{n}}.
\end{eqnarray*}

%

Moreover consider  ${\cal D}_n$ as the following event:
$$
{\cal D}_n= \left\{ D_{n} \leq \frac{e^{-n\epsilon^{2}_{n}}}{2} \Pi_{n} \left (  {\cal U}(q_{0},\epsilon_{n}) \right) \right\}.$$
By means of the test procedure defined in \eqref{test}, $\psi(\mathcal{H}_n)$, $\mathbb{E}_{0}^{(n)}(A_{1})$ could be written as follows
\begin{eqnarray}
\nonumber
\mathbb{E}_{0}^{(n)}(A_{1})&=&\mathbb{E}_{0}^{(n)}\Bigl{(}\frac{N_{n}}{D_{n}}\Bigr{)}\\
\nonumber
&\leq &\mathbb{E}_{0}^{(n)}[\psi(\mathcal{H}_n)]+\mathbb{E}_{0}^{(n)}\Big{[}(1-\psi(\mathcal{H}_n))\frac{N_{n}}{D_{n}} \Big{\{}\BBone_{ {\cal D}_n}
+\BBone_{{\cal D}_n^{\complement}}\Big{\}} \Big{]} \\
\nonumber
&\leq& \mathbb{E}_{0}^{(n)}[\psi(\mathcal{H}_n)]+\mathbb{E}_{0}^{(n)}\Big{[}(1-\psi(\mathcal{H}_n))\frac{N_{n}}{D_{n}}\BBone_{{\cal D}_n^{\complement}}\Big{]} 
+ \mathbb{P}^{(n)}_{0}\Big{(}{\cal D}_n\Big{)}\\
&:=&T_{1}+T_{2}+T_{3}. \label{eq:T1T2T3}
\end{eqnarray}

\vspace{0.2cm}

To bound from above $\mathbb{E}_{0}^{(n)}(A_{1})$, it is sufficient to upper bound
every term in the right-hand side of (\ref{eq:T1T2T3}).

\begin{itemize}
\item \textsc{Term $T_1$.}  We apply Corollary \ref{cor:test-global} and obtain that 
there exists $K>0$ such that 
\begin{align}
T_{1}&=\mathbb{E}^{(n)}_{0}[\psi(\mathcal{H}_{n})]\leq e^{-K n\epsilon_n^2M^2}. \label{T1} 
\end{align}

\item \textsc{Term $T_2$.}
We apply once again Corollary \ref{cor:test-global}, which combined with \textbf{H3} entails that there exists ${\tilde K}>0$ such that
\begin{align}
T_{2} 
&\leq\int_{B_{d_{\nu^*}}^{\complement}(q_{0},\epsilon_{n}M)\cap \mathcal{Q}_n}\mathbb{E}_q^{(n)}[1-\psi(\mathcal{H}_{n})]d\Pi_{n}(q)\frac{2}{e^{-n\epsilon_{n}^{2}}{\Pi_{n}}\bigl{(}{\cal U}(q_{0},\epsilon_{n})\bigr{)}} \nonumber\\
&\leq  \sup_{q  \in B_{d_{\nu^*}}^{\complement}(q_{0},\epsilon_{n}M)\cap \mathcal{Q}_n}{\mathbb{E}_q^{(n)}[1-\psi(\mathcal{H}_{n})]\frac{2}{e^{-n\epsilon_{n}^{2}}{\Pi_{n}}\bigl{(}{\cal U}(q_{0},\epsilon_{n})\bigr{)}}} \nonumber \\
&\leq  e^{-{\tilde K} n \epsilon_n^2 M^{2}} \frac{2}{e^{-n\epsilon_{n}^{2}}{\Pi_{n}}\bigl{(}{\cal U}(q_{0},\epsilon_{n})\bigr{)}} \nonumber\\
&\leq 2e^{-({\tilde K} M^2 -1-c)n \epsilon_n^2  } \leq 2e^{- \kappa n \epsilon_n^2  }, \label{T2}
\end{align}
where $\kappa:={\tilde K} M^2 -1-c$ is positive under the condition that $M$ is sufficiently large.


\item \textsc{Term $T_3$.}
Consider the following subspace of ${\cal Q}$ 
\begin{equation*}
\mathcal{V}_{n}:=\Big{\{}q \in \mathcal{Q}:
\log\frac{{\cal L}_n(q)}{{\cal L}_n(q_0)}+K(\mathbb{P}^{(n)}_{0},\mathbb{P}_q^{(n)})\geq\frac{n\epsilon^{2}_{n}}{2}\Big{\}},
\end{equation*}
%
and observe that 
\begin{eqnarray*}
D_{n}
&\geq &\int_{{\cal U}(q_{0},\epsilon_{n})\cap{\mathcal{V}_{n}}}^{}\exp\Bigl{(}\log\frac{{\cal L}_n(q)}{{\cal L}_n(q_0)} + K(\mathbb{P}^{(n)}_{0},\mathbb{P}_q^{(n)}) -K(\mathbb{P}^{(n)}_{0},\mathbb{P}_q^{(n)})\Bigr{)} d\Pi_{n}(q)\\
&\geq&\exp\Bigl{(}\frac{-n\epsilon^{2}_{n}}{2}\Bigr{)}{\Pi_{n}}\Bigl{(}{\cal U}(q_{0},\epsilon_{n})\cap{\mathcal{V}_{n}}\Bigr{)}.
\end{eqnarray*}



\vspace{0.2cm}

It then follows from Fubini's theorem and Markov's inequality that
\begin{align}
T_{3}&\leq
\mathbb{P}_{0}^{(n)}\Big{(}e^{\frac{-n\epsilon^{2}_{n}}{2}}{\Pi_{n}}\Bigl{(}{\cal U}(q_{0},\epsilon_{n})\cap{\mathcal{V}_{n}}\Bigr{)}\leq{{\frac{e^{-n\epsilon_{n}^{2}}}{2}}{{\Pi_{n}}\bigl{(}{\cal U}(q_{0},\epsilon_{n})\bigr{)}}}\Bigr{)} \nonumber\\
&={\mathbb{P}_{0}^{(n)}\Big{(}{\Pi_{n}}\Bigl{(}{\cal U}(q_{0},\epsilon_{n})\cap{\mathcal{V}^{\complement}_{n}}\Bigr{)}\geq{\bigl{(}1-\frac{1}{2}e^{\frac{-n\epsilon_{n}^{2}}{2}}\bigr{)}{{\Pi_{n}}\bigl{(}{\cal U}(q_{0},\epsilon_{n})\bigr{)}\Bigr{)}}}} \nonumber \\
&\leq\frac{2}{\bigl{(}2-e^{\frac{-n\epsilon_{n}^2}{2}}\bigr{)}{{\Pi_{n}}\bigl{(}{\cal U}(q_{0},\epsilon_{n})\bigr{)}}}\mathbb{E}^{(n)}_{0}\Bigl{(}\Pi_n\bigl{(}\mathcal{V}^{\complement}_{n}\cap{{\cal U}}(q_{0},\epsilon_{n})\bigr{)}\Bigr{)} \nonumber \\
&\leq \frac{2}{\bigl{(}2-e^{\frac{-n\epsilon_{n}^2}{2}}\bigr{)}{{\Pi_{n}}\bigl{(}{\cal U}(q_{0},\epsilon_{n})\bigr{)}}}
\times\int_{{\cal U}(q_{0},\epsilon_{n})}{\mathbb{P}^{(n)}_{0}\Bigl{(}\vert \log\frac{{\cal L}_n(q_0)}{{\cal L}_n(q)} -K(\mathbb{P}^{(n)}_{0},\mathbb{P}^{(n)}_{q})\vert>\frac{n\epsilon^{2}_{n}}{2}}\Bigr{)}d\Pi_{n}(q) \nonumber \\
&\leq\frac{2}{\bigl{(}2-e^{\frac{-n\epsilon_{n}^2}{2}}\bigr{)}{{\Pi_{n}}\bigl{(}{\cal U}(q_{0},\epsilon_{n})\bigr{)}}}
\int_{{\cal U}(q_{0},\epsilon_{n})}^{}{\mathbb{V}_0(\mathbb{P}^{(n)}_{0},\mathbb{P}_q^{(n)})d\Pi_{n}(q)\frac{4}{n^{2}\epsilon^{4}_{n}}} \nonumber \\
&\leq \frac{8}{n\epsilon^{2}_{n}\bigl{(}2-e^{\frac{-n\epsilon_{n}^2}{2}}\bigr{)}}.\label{T3}
\end{align}


\end{itemize}

\bigskip
Third, let us  turn  to $A_2$ which is rewritten as follows
\begin{eqnarray*}
A_{2}&=&
\frac{\int_{B_{d_{\nu^*}}^{\complement}(q_{0},\epsilon_{n}M)\cap{\mathcal{Q}_n^{\complement}}} \frac{{\cal L}_n(q)}{{\cal L}_n(q_0)}d{\Pi_{n}}(q)}
{\int_{\mathcal{Q}}^{}\frac{{\cal L}_n(q)}{{\cal L}_n(q_0)}d{\Pi_{n}}(q)} :=\frac{\widetilde{N}_{n}}{D_{n}}.\\
\end{eqnarray*}


Then, using  Equation \eqref{T3} and from Assumptions \textbf{H3} and \textbf{H4}, we obtain
\begin{align}
\mathbb{E}_{0}^{(n)}(A_{2}) 
&=
\mathbb{E}_{0}^{(n)}      \Big{(}
\frac{\widetilde{N}_{n}}{D_{n}}
\Big{\{}
\BBone_{D_{n}\leq \frac{e^{-n\epsilon_n^2}}{2}
\Pi_{n} \bigl{(} {\cal U}(q_{0},\epsilon_{n})\bigr{)}}+\BBone_{D_{n}> \frac{e^{-n\epsilon_n^2}}{2}
\Pi_{n} \bigl{(}{\cal U}(q_{0},\epsilon_{n}) \bigr{)}}
\Big{\}} \Big{)}  \nonumber \\
&\leq \mathbb{P}_{0}^{(n)}
\Bigl{(}   {\cal D}_n \Bigr{)}+
\mathbb{E}^{(n)}_{0} \big{(} \widetilde{N}_{n} \big{)}  \frac{2}{e^{-n\epsilon^{2}_{n}}\Pi_{n} \bigl{(} {\cal U}(q_{0}, \epsilon_{n})\bigr{)}} \nonumber \\
&\leq \mathbb{P}_{0}^{(n)}
\Bigl{(}   {\cal D}_n \Bigr{)}+ \Pi_{n}\bigl{(}\mathcal{Q}_n^{\complement}\bigr{)}\frac{2}{e^{-n\epsilon^{2}_{n}}\Pi_{n} \bigl{(} {\cal U}(q_{0}, \epsilon_{n})\bigr{)}} \nonumber \\
& \leq \frac{8}{n\epsilon^{2}_{n}\bigl{(}2-e^{\frac{-n\epsilon_{n}^2}{2}}\bigr{)}} + 2e^{-(c+1) n \epsilon_n^2}. \label{A2}
\end{align}

Finally, Inequalities \eqref{T1}--\eqref{A2} lead to the desired result \eqref{eq:resu}.

\vspace{0.5cm}

\smallskip
\noindent \textbf{Acknowledgments}

The authors would like to thank the F\'{e}d\'{e}ration de Recherche Math\'{e}matiques des Pays de Loire
(FR 2962) for their support.



%
%
%

\end{document}